\documentstyle[amssymb,amsfonts]{amsart}

 at 10 true pt

\def\be#1{ \begin{equation}\label{#1} }

\def\bas{\begin{align*}}
\def\eas{\end{align*}}
\def\bi{\begin{itemize}}
\def\ei{\end{itemize}}

\newenvironment{proof}{\noindent {\bf Proof} }{\endprf\par}
\def \endprf{\hfill  {\vrule height6pt width6pt depth0pt}\medskip}
\def\emph#1{{\it #1}}
\def\textbf#1{{\bf #1}}

\def\BZ{{\mathbf Z}}

\def\ep{{\epsilon}}

\theoremstyle{plain}
  \newtheorem{theorem}[subsection]{Theorem}
  \newtheorem{question}[subsection]{Question}

  \newtheorem{lemma}[subsection]{Lemma}

\theoremstyle{remark}

\theoremstyle{definition}

\include{psfig}

\begin{document}

\title[Subset sums in $\BZ_p$]
{Subset sums in $\BZ_p$}

\author{Hoi H. Nguyen}
\address{Department of Mathematics, Rutgers, Piscataway, NJ 08854}
\email{}
\thanks{}

\author{Endre Szemer\'edi}
\address{Department of Computer Science, Rutgers, Piscataway, NJ 08854}
\email{}
\thanks{}

\author{Van H. Vu}
\address{Department of Mathematics, Rutgers, Piscataway, NJ 08854}
\email{vanvu@@math.rutgers.edu}
\thanks{V. Vu is an A. Sloan Fellow and is supported by an NSF Career Grant.}

\begin{abstract}
Let $\BZ_p$ be the finite field of prime order $p$ and $A$ be a
subset of $\BZ_p$. We prove several sharp results about the
following two basic questions:

(1) When can one represent zero as a sum of distinct elements of
$A$ ?

(2) When can one represent every element of $\BZ_p$ as a sum  of
distinct elements of $A$ ?

\end{abstract}

\maketitle

\section{Introduction}
Let $A$ be an additive group and $A$ be a subset of $A$. We denote
by $S_A$ the collection of subset sums of $A$:

$$
S_A= \{ \sum_{x\in B}x |B\subset A, |B| < \infty \}.
$$

The following two questions are among the most popular questions
in additive combinatorics

\begin{question} \label{question:1} When  $ 0 \in S_A$  ?
\end{question}

\begin{question} \label{question:2} When  $S_A= G$ ?
\end{question}

If $S_A$ does not contain the zero element, we say that $A$ is
{\it zero-sum-free}. If $S_A =G$ ($S_A \neq G$), then we say that
$A$ is {\it complete (incomplete)}.

In this paper, we focus on the case $G=\BZ_p$, the cyclic group of
order $p$, where $p$ is a large prime. The asymptotic notation
will be used under the assumption that $p \rightarrow \infty$. For
$x \in \BZ_p$, $\|x\|$ (the norm of $x$) is the distance from $x$
to $0$. (For example, the norm of $p-1$ is 1.) All logarithms have
natural base and $[a,b]$ denotes the set of integers between $a$
and $b$.

\subsection{A sharp bound on the maximum cardinality of a zero-sum-free set}

How big can a zero-sum-free set be ? This question was raised by
Erd\H os and Heilbronn \cite{EH} in 1964. In \cite{Szem1},
Szemer\'edi proved that

\begin{theorem} \label{theorem:Szem} There is a positive constant
$c$ such that the following holds. If $A \subset \BZ_p$ and $|A|
\ge cp^{1/2}$, then $0 \in S_A$. \end{theorem}

A result of Olson \cite{Olson} implies that one can set $c=2$. More
than a quarter of century  later, Hamindoune and Z\'emor \cite{HZ}
showed that one can set $c= \sqrt 2 +o(1)$, which is asymptotically
tight.

\begin{theorem} \label{theorem:HZ}  If $A \subset \BZ_p$ and $|A|
\ge (2p)^{1/2} +5 \log p$, then $0 \in S_A$. \end{theorem}

Our first result removes  the logarithmic term in Theorem
\ref{theorem:HZ}, giving  the best possible bound (for all
sufficiently large $p$). Let $n(p)$ denote the largest integer such
that $\sum_{i=1}^{n-1}i < p$.

\begin{theorem} \label{theorem:main1} There is a constant $C$ such
that the following holds for all prime $p \ge C$.

\begin{itemize}

\item  If $p \neq \frac{n(p)(n(p) +1)}{2} -1$, and $A$ is a subset
of $\BZ_p$ with $n(p)$ elements, then $0 \in S_A$.

\item If $p = \frac{n(p)(n(p) +1)}{2} -1$, and $A$ is a subset of
$\BZ_p$ with $n(p)+1$ elements, then $0 \in S_A$. Furthermore, up
to a dilation, the only $0$-sum-free set with $n(p)$ elements is
$\{-2,1,3,4, \dots, n(p)\}$.

\end{itemize}

\end{theorem}
To see that the bound in the first case is sharp, consider $A=\{1,
2, \dots, n(p)-1\}$.

\subsection{The structure of zero-sum-free sets with cardinality closed to maximum}

 Theorem \ref{theorem:main1} does not provide information about
zero-sum-free sets of size slightly smaller than $n(p)$. The arch
typical example for a zero-sum-free set is a set whose sum of
elements (as positive integers between 1 and $p-1$) is less than
$p$. The general phenomenon we would like to support here is that
a zero-sum-free set with sufficiently large cardinality should be
close to such a set. In \cite{D1}, Deshouillers \cite{D1} showed

\begin{theorem} \label{theorem:D1} Let $A$ be a zero-sum-free
subset of $\BZ_p$  of size at least $p^{1/2}$. Then there is some
non-zero element $b \in \BZ_p$ such that
$$\sum_{a \in bA, a<p/2} \| a\| \le p+O(p^{3/4}\log p) $$
and
$$\sum_{a \in bA, a>p/2} \| a\|=O(p^{3/4} \log p).$$ \end{theorem}

The main issue here is the magnitude of the error term.  In the
same paper, there is a construction of a zero-sum-free set with
$cp^{1/2}$ elements ($c >1)$ where

$$\sum_{a \in bA, a<p/2} \| a\| = p+\Omega(p^{1/2}) $$
and
$$\sum_{a \in bA, a>p/2} \| a\|= \Omega (p^{1/2}).$$

\noindent It is conjectured \cite{D1}  that  $p^{1/2}$ is the
right order of magnitude of the error term. Here we confirm this
conjecture, assuming that $|A|$ is sufficiently close to the upper
bound.

\begin{theorem}\label{theorem:main2} Let $A$ be a zero-sum-free
subset of $\BZ_p$  of size at least $.99 (2p)^{1/2}$. Then there is
some non-zero element $b \in \BZ_p$ such that
$$\sum_{a \in bA, a<p/2} \| a\| \le p+O(p^{1/2}) $$
and
$$\sum_{a \in bA, a>p/2} \| a\|=O(p^{1/2}).$$
\end{theorem}

The constant $.99$ is adhoc and can be improved. However, we do
not elaborate on this point.

\subsection{Complete sets}

All questions concerning zero-sum-free sets are also natural for
incomplete sets.  Here is a well-known result of Olson
\cite{Olson}

\begin{theorem} Let $A$ be a subset of $\BZ_p$ of more than
$(4p-3)^{1/2}$ elements, then $A$ is complete. \end{theorem}

Olson's bound is essentially sharp. To see this, observe that if
the sum of the norms of the elements of $A$ is less than $p$, then
$A$ is incomplete.  Let $m(p)$ be the largest cardinality of a
small set. One can easily verify  that $m(p)= 2p^{1/2} +O(1)$. We
now want to study the structure of incomplete sets of size close
to $2p^{1/2}$.  Deshouillers and Freiman \cite{DF} proved

\begin{theorem} \label{theorem:DF} Let $A$ be an incomplete
subset of $\BZ_p$  of size at least $(2p)^{1/2}$. Then there is
some non-zero element $b \in \BZ_p$ such that

$$\sum_{a \in b  A} \| a\| \le p +O(p^{3/4} \log p). $$

\end{theorem}

\noindent Similar to the situation with Theorem \ref{theorem:D1},
it is conjectured that the right error term has order $p^{1/2}$
(see \cite{D2} for a construction that matches this bound from
below). We establish this conjecture for sufficiently large $A$.

\begin{theorem} \label{theorem:main3} Let $A$ be an incomplete
subset of $\BZ_p$  of size at least $1.99 p^{1/2}$. Then there is
some non-zero element $b \in \BZ_p$ such that

$$\sum_{a \in b  A} \| a\| \le p +O(p^{1/2}). $$

\end{theorem}

{\it Added in proof.} While this paper was written, Deshouillers
informed us that he and Prakash have obtained a result similar to
Theorem \ref{theorem:main1}.

\section{Main lemmas}

The main tools in our proofs are the following results from
\cite{SzemVu1}.

\begin{theorem} \label{lemma:main1} Let  $A$ be a zero-free-sum subset of
$\BZ_p$. Then we can partition $A $ into two disjoint sets $ A'$
and $ A^{''} $ where

\begin{itemize}

\item $A'$ has negligible  cardinality: $|A'| = O(p^{1/2} /\log^2
p).$

\item The sum of the elements of (a dilate of) $A^{''}$ is small:
There is a non-zero element $b \in \BZ_p$ such that the elements
of $b  A^{''} $ belong to the interval $[1, (p-1)/2]$ and their
sum is less than $p$.

\end{itemize}
\end{theorem}

\begin{theorem}\label{lemma:main2}
Let  $A$ be an incomplete subset of $\BZ_p$. Then we can partition
$A $ into two disjoint sets $ A'$ and $ A^{''} $ where

\begin{itemize}

\item $A'$ has negligible  cardinality: $|A'| = O(p^{1/2} /\log^2
p).
$

\item The norm sum of the elements of (a dilate of) $A^{''}$ is
small: There is a non-zero element $b \in \BZ_p$ such that the sum
of the norms of the  elements of $b   A^{''} $ is less than $p$.

\end{itemize}
\end{theorem}

The above two theorems were proved (without being formally stated)
in \cite{SZemVu1}. A stronger version of these theorems will
appear in a forth coming paper \cite{NgSzV}. We also need the
following simple lemmas.

\begin{lemma}\label{lemma:simple4} Let $T'\subset T$ be sets of
integers with the following property. There are integers $a \le b$
such that $[a,b] \subset S_{T'}$ and the non-negative (non-positive)
elements of $T\backslash T'$ are less than $b-a$ (greater than
$a-b$). Then

$$[a, b+ \sum_{x \in T\backslash T', x \ge 0} x] \subset S_T. $$

$$([a+ \sum_{x \in T\backslash T', x \le 0} x,b] \subset S_T. )$$

\end{lemma}

\noindent The (almost trivial) proof is left as an exercise.

\begin{lemma} \label{lemma:simple5} Let $K=\{k_1, \dots, k_l\}$ be
a subset of $\BZ_p$, where the $k_i$ are positive integers and
$\sum_{i=1}^l k_i \le  p$. Then $|S_K| \ge l(l+1)/2.$

\end{lemma}

To verify this lemma, notice that the numbers

$$k_1, \dots, k_l, k_1+k_l, k_2+k_l, \dots, k_{l-1}+k_l, k_1 +
k_{l-1} +k_l, \dots, k_{l-2} + k_{l-1}+k_l, \dots, k_1+\dots + k_l
$$

\noindent are different and all belong to $S_K$.

\section{Proof of Theorem \ref{theorem:main1}}

Let $A$ be a zero-free-sum subset of $\BZ_p$ with size $n(p)$. In
fact, as there is no danger for misunderstanding, we will write $n$
instead of $n(p)$. We start with few simple observations.

Consider the partition $A= A' \cup A^{''}$ provided by Theorem
\ref{lemma:main1}. Without loss of generality, we can assume that
the element $b$ equals one. Thus $A^{''} \subset [1, (p-1)/2]$ and
the sum of its elements is less than $p$. Set $I_n:=[1,n]$ be the
set of the first $n$ positive integers. We first show that most of
the elements of $A^{''}$ belong to $I_n$.

\begin{lemma} \label{lemma:simple1} $|A^{''} \cap I_n | \ge n -
O(n / \log n). $ \end{lemma}

\begin{proof} By the definition of $n$ and the property of
$A^{''}$

$$\sum_{i=1}^n i \ge p > \sum_{a \in A^{''}} a. $$

\noindent Assume that $A^{''}$ has $l$ elements in $I_n$ and $k$
elements outside. Then

$$  \sum_{a \in A^{''}}
 a \ge \sum_{i=1}^l i + \sum_{j=1}^k (n+j).
$$

\noindent It follows that

$$\sum_{i=1}^n i  >  \sum_{i=1}^l i + \sum_{j=1}^k (n+j), $$

\noindent which, after a routine simplification, yields

$$ (l+n+1) (n-l)  > (2n+k) k. $$

\noindent On the other hand, $n \ge k+l= |A^{''}| \ge
n-O(n/\log^2n)$, thus $n-l =k+O(n/\log^2 n)$ and $n+l+1 \le 2n-k+1$.
So there is a constant $c$ such that

$$(2n-k+1) (k+ cn/\log^2 n) > (2n+k) k, $$

\noindent or equivalently

$$\frac{cn}{k \log^2 n } > \frac{k+1}{2n-k+1}. $$

\noindent Since $2n-k+1\le 2n+1$, a routine consideration shows that
$k^2 \log^2 n =O(n^2)$ and thus $k= O(n/\log n)$, completing the
proof.
\end{proof}

The above lemma shows that most of the elements of $A^{''}$ (and
$A$) belong to $I_n$. Let $A_1 = A \cap I_n$. It is trivial that

$$|A_1| \ge |A^{''} \cap I_n| = n -O(n /\log n). $$

\noindent Let $A_2 = A \backslash A_1$. We have

$$t := |I_n \backslash A_1| = |A_2| =|A| -|A_1| = O(n /\log n). $$

\noindent Next we show that $S_{A_1}$ contains a very long interval.
Set $I:= [2t+3, (n+1) (\lfloor n/2 \rfloor -t-1)]$. The length of
$I$ is $(1-o(1))p$; thus $I$ almost cover $\BZ_p$.

\begin{lemma} \label{lemma:simple2} $I \subset S_{A_1}.$
\end{lemma}

\begin{proof}  We
need to show that every element $x$ of in this interval can be
written as a sum of distinct elements of $A_1$. There are two
cases:

{\bf Case 1.} $2t+3 \le x \le n.$ In this case $A_1$ contains  at
least $x-1-t \ge (x+1)/2$ elements in the interval $[1,x-1]$. This
guarantees that there are two distinct elements of $A_1$ adding up
to $x$.

{\bf Case 2.} $x= k(n+1) +r$ for some $1 \le k \le \lfloor n/2
\rfloor -t-2$  and $0 \le r \le n+1$. First, notice that since
$|A_1|$ is very close to $n$ (in fact it is enough to have $|A_1|$
slightly larger than $2n/3$ here), one can find three distinct
elements $a,b,c \in A_1$ such that $a+b+c = n+1+r$. Consider the set
$A_1'= A_1\backslash \{a,b,c\}$. We will represent $x-(n+1+r) =(k-1)
(n+1)$ as a sum of distinct elements of $A_1'$.
 Notice that there are exactly $\lfloor n/2 \rfloor$ ways to
write $n+1$ as a sum of two different positive integers. We
discard a pair if (at least) one of its two elements is not in
$A_1'$. Since $|A'_1|= n-t-3$, we discard at most $t+3$ pairs. So
there are at least $\lfloor n/2 \rfloor -t-3$ different pairs
$(a_i,b_i)$ where $a_i, b_i \in A_1'$ and $a_i+b_i =(n+1)$. Thus,
$(k-1)(n+1)$ can be written as  a sum of distinct pairs. Finally,
$x$ can be written as a sum of $a,b,c$ with these pairs.
\end{proof}

Now we investigate the set $A_2= A\backslash A_1$. This is the
collection of elements of $A$ outside the interval $I_n$. Since
$A$ is zero sum free, $0 \notin A_2 +I $ thanks to Lemma
\ref{lemma:simple2}. It follows that

$$A_2 \subset \BZ_p \backslash \{I_n \cup (-I) \cup \{0\} \subset J_1 \cup J_2, $$

\noindent where $J_1:= [-2t-2,-1]$ and $J_2= [(n+1), p
-(n+1)(\lfloor n/2 \rfloor -t-t)] =[(n+1),q]$. We set $B:=A_2 \cap
J_1$ and $C:=A_2 \cap J_2$.

\begin{lemma}\label{lemma:simple3} $S_B \subset J_1.$ \end{lemma}

\begin{proof} Assume otherwise. Then there is a subset $B'$ of $B$
such that $\sum_{a \in B'} a \le -2t-3$ (here the elements of $B$
are viewed as negative integers between $-1$ and $-2t-3$). Among
such $B'$, take one where $\sum_{a \in B'} a$ has the smallest
absolute value. For this $B'$, $\le -4t-4 \sum_{a \in B'} a \le
-2t-3$. On the other hand, by Lemma \ref{lemma:simple2}, the
interval $2t+3, 4t+4$ belongs to $S_{A_1}$. This implies that $0 \in
S_{A_1} + S_b \subset S_A$, a contradiction. \end{proof}

Lemma \ref{lemma:simple3} implies that $\sum_{a \in B} |a| \le
2t+2$, which yields

\begin{equation} \label{boundonB} |B| \le 2 (t+1)^{1/2}.
\end{equation}

Set $s:=|C|$. We have $s \ge t - 2 (t+1)^{1/2}. $Let $c_1 < \dots
< c_s$ be the elements of $C$ and $h_1 < \dots < h_t$ be the
elements of $I_n \backslash A_1$.

By the definition of $n$, $\sum_{i=1}^n i > p > \sum_{i=1}^{n-1}
i$. Thus, there is an (unique) $h \in I_n$ such that

\begin{equation} \label{missingh} p= 1+\dots + (h-1) +(h+1) + \dots + n.
\end{equation}

\noindent A quantity which plays an important  role in what
follows is

$$ D:= \sum_{i=1}^s c_i -\sum_{j=1}^{t} h_j. $$

\noindent Notice that if we replace the $h_j$ by the $c_i$ in
\eqref{missingh}, we represent $p+D$ as a sum of distinct elements
of $A$

\begin{equation} \label{missingh1} p+D = \sum_{a\in X, X \subset
A}  a. \end{equation}

\noindent The leading idea now is to try to cancel $D$ by throwing
a few elements  from the right hand side or adding a few negative
elements (of $A$) or both. If this was  always possible, then we
would have  a representation of $p$ as a sum of distinct elements
in $A$ (in other words $0 \in S_A$), a contradiction. To conclude
the proof of Theorem \ref{theorem:main1}, we are going to show
that the only case when it is not possible is when $p =
n(n+1)/2-1$ and $A=\{-2,1,3,4, \dots, n \}$. We consider two
cases:

{\bf Case 1. $h \in A_1$.}  Set $A_1'= A_1 \backslash \{h \}$ and
apply Lemma \ref{lemma:simple2} to $A_1'$, we conclude that
$S_{A_1'}$ contains the  interval $I'=[2(t+1)+3, (n+1)(\lfloor n/2
\rfloor -t-2)]$.

\noindent

\begin{lemma}\label{lemma:simple6} $D < 2(t+1)+3. $
\end{lemma}

\begin{proof} Assume $D \ge 2(t+1)+3$. Notice that the largest
element in $J_2$ (and thus in $C$) is less than the length of
$I'$. So by removing the $c_i$ one by one from $D$, one can obtain
a sum $D'= \sum_{i=1}^{s'} c_i - \sum_{j=1}^{t} h_j$ which belongs
to $I'$, for some $s' \le s$. This implies

$$\sum_{i=1}^{s'} c_i = \sum_{j=1}^t h_j + \sum_{a \in X} a $$

\noindent for some subset $X$ of $A_1'$. Since $h \notin A_1'$, the
right hand side is a sub sum of the right hand side  of
\eqref{missingh}. Let $Y$ be the collection of the missing elements
(from the right hand side  of \eqref{missingh}). Then $Y \subset
A_1$ and $\sum_{i=1}^{s'} c_i + \sum_{a \in Y} a = p$. On the other
hand, the left hand side belongs to $S_{A_1} + S_{A_2} \subset S_A$.
It follows that $0 \in S_A$, a contradiction.
\end{proof}

\noindent Now we take a close look at the inequality $D < 2(t+1)+
3$. First, observe that since $A$ is zero-sum-free, $-S_B \subset
\{h_1, \dots, h_t \}$. By Lemma \ref{lemma:simple3}, $\sum_{a \in
B} |a| \le 2t+2 <p$. As $B$ has $t-s$ elements, by Lemma
\ref{lemma:simple5}, $S_B$ has at least $(t-s)(t-s+1)/2$ elements.
It follows that

$$\sum_{i=1}^t h_i \le (2t+2) + \sum_{j=0}^ {(t-(t-s)(t-s+1)/2) +1}(n- j). $$

\noindent On the other hand, as all elements of $C$ are larger
than $n$

$$\sum_{i=1}^s c_s \ge \sum_{i=1}^s (n+i). $$

\noindent It follows that $D$ is at least

$$ \sum_{i=1}^s (n+i) -  (2t+2) - \sum_{j=0} ^{ (t-(t-s)(t-s+1)/2) +1} (n-j) . $$

\noindent If $t-s \ge 2$, then $s > t- (t-s)(t-s+1)/2$, so the
last formula has order $\Omega (n) \gg t$, thus $D \gg 2(t+1)+3$,
a contradiction. Therefore, $t-s$ is either $0$ or $1$.

If $t-s=0$, then $D =\sum_{i=1}^t c_i -\sum_{i=1}^t h_i \ge t^2$.
This is larger than $2t+5$ if $t \ge 4$. Thus, we have
$t=0,1,2,3$.

\begin{itemize}

\item $t=0.$ In this case $A=I_n$ and $0 \in S_A$.

\item $t=1$. In this case $A= I_n \backslash \{h_1 \} \cup c_1$.
If $c_1 -h_1 \neq h$, then we could substitute $c_1$ for $h_1 +
(c_1-h_1)$ in \eqref{missingh} and have $0 \in S_A$. This means
that $h=c_1-h_1$. Furthermore, $h < 2t+5=7$ so both $c_1$ and
$h_1$ are close to $n$.  If $h \ge 3$,

$$p= \sum_{i=1}^{h-1} i +\sum_{j=h+1}^n j = \sum_{i=2}^{h-2} i +
\sum_{h+1 \le j \le n, j \neq h_1} j + c_ 1.$$  Similarly, if
$h=1$ or $2$  then we have

$$p =\sum_{i=1}^h i + \sum_{h+2 \le j \le n, j \neq h_1} j + c_1. $$

\item $t > 1$. Since $D < 2t +5$, $h_1, \dots, h_t$ are all larger
than $n-2t-4$. As $p$ is sufficiently large, we can assume $n \ge
4t +10$, which implies that  $[1, 2t+5] \subset A_1$.  If $h \neq
1$, then it is easy to see that $[3, 2t+5] \subset S_{A_1
\backslash \{h\}}$. As $t >1$, $D \ge t^2 \ge  4$ and can be
represented as a sum of elements in ${A_1 \backslash \{h\}}$.
Omitting these elements from \eqref{missingh1}, we obtain a
representation of $p$ as a sum of elements of $A$. The only case
left is $h=1$ and $D=4$. But $D$ can equal $4$ if and only if
$t=2$, $c_1=n+1, c_2 = n+2, h_1=n-1, h_2=n$. In this case, we have

$$p= \sum_{i=2}^n i = 2+3 + \sum_{i=5}^{n+2} i .$$

\end{itemize}

Now we turn to the case $t-s=1$. In this case $B$ has exactly one
element  in the interval $[-2t-2,-1]$ (modulo $p$) and $D$ is at
least $s^2-(2t+2)= (t-1)^2 -(2t+2)$. Since $D < 2t +5$, we
conclude that $t$ is at most 6. Let $-b$ be the element in $B$
(where $b$ is a positive integer). We have $b \le 2t+2 \le 14$.
$A_1$ misses exactly $t$ elements from $I_n$; one of them is $b$
and all other are close to $n$ (at least $n-(2t+4)$). Using this
information, we can reduce the bound on $b$ further. Notice that
the whole interval $[1,b-1]$ belongs to $A_1$. So if $b \ge 3$,
then there are two elements $x,y$ of $A_1$ such that $x+y=b$. Then
$x+y+(-b)=0$, meaning $0 \in S_A$. It thus remain to consider
$b=1$ or $2$. Now we consider a few cases depending on the value
of $D$. Notice that $D \ge s^2-b \ge -2$. In fact, if $s \ge 2$
then $D \ge 2$. Furthermore, if $s=0$, then $t=1$ and $D= -h_1
=-b$.

\begin{itemize}

\item $D \ge 5$. Since $A_1$ misses at most one element in $[1,D]$
(the possible missing element is $b$), there are two elements of
$A_1$ adding up to $D$. Omitting these elements from
\eqref{missingh1}, we obtain a representation of $p$ as a sum of
distinct elements of $A$.

\item $D=4$. If $b =1$, write $p= \sum_{a \in X, a \neq 2 } a +
(-b)$. If $b=2$, then $p= \sum_{a \in X, a \neq 1,3} a $. (Here
and later $X$ is the set in \eqref{missingh1}.)

\item $D = 3$. Write $p= \sum_{a \in X, a \neq 3-b} a + (-b)$.

\item $D=2$.  If $b=1$, then $p= \sum_{a \in X, a \neq 2} a$. If
$b=2$, then $p= \sum_{a\in X} a +(-2). $

\item $D=1$. If $b=1$, then $p= \sum_{a \in X} a +(-1)$. If $b=2$,
then $p= \sum_{a\in X, a \neq 1} a. $

\item $D=0$. In this case \eqref{missingh1} already provides a
representation of $p$.

\item $D=-1$. In this case $s < 2$. But since $h \neq b$, $s$
cannot be $0$. If $s=1$ then $b=2$ and $c_1=n+1$, $h_1= n$. By
\eqref{missingh}, we have $p= \sum_{i=1}^{h-1} i + \sum_{j=h+1}^n
j$ and so

$$p +(h-1) = \sum_{1 \le i \le n+1, i \notin \{2, n\}} i $$ where
the right hand side consists of elements of $A$ only. If $h-1 \in
A$ then we simply omit it from the sum. If $h-1 \notin A$, then
$h-1=2$ and $h=3$. In this case, we can write

$$p= \sum_{1 \le i \le n+1, i \notin \{2, n\}} i +(-2).$$

\item $D=-2$. This could only occur if $s=0$ and $b=2$. In this
case $A= \{-2, 1,3, \dots, n \}$. If $h=1$, then $p=\sum_{i=2}^n =
n(n+1)/2 -1$ and we end up with the only exceptional set. If $h \ge
3$, then $p+(h-2) =\sum_{1 \le i \le n, i \neq 2} i$. If $ h\neq 4$,
then we can omit $h-2$ from the right hand side to obtain a
representation of $p$. If $h=4$, then we can write

$$p= \sum_{1 \le i \le n, i \neq 2} i + (-2). $$

\end{itemize}

\vskip2mm

{\bf Case 2. $h \notin A$}. In this case we can consider $A_1$
instead of $A_1'$. The consideration is similar and actually
simpler. Since $h \notin A$, we only need to consider $D:=
\sum_{i=1}^s c_i - \sum_{1 \le j \le t, h_j \neq h} h_j$.
Furthermore, as $h \notin A$, if $s=0$ we should have $h=b$ and
this forbid us to have any exceptional structure in the case
$D=-2$. The detail is left as an exercise.

\section{Proof of Theorem \ref{theorem:main2}}

\noindent We follow the same terminology used in the previous
section. Assume that $A$ is zero-sum-free and $|A|=\lambda n=\lambda
(2p)^{1/2}$ with some $1\ge \lambda\ge .99.$ Furthermore, assume
that the element $b$ in Theorem \ref{lemma:main1} is one. We will
use the notation of the previous proof. Let the {\it core} of $A$ be
the collection of all pairs $(a,a') \in A \times A, a \neq a'$ and
$a+a'=n+1$. Theorem \ref{theorem:main2} follows directly from the
following two lemmas.

\begin{lemma} \label{lemma:main2-1} The core of $A$ has size at
least $.6n$.
\end{lemma}

\begin{lemma} \label{lemma:main2-2} Let $A$ be a zero-sum-free
set whose core has size at least $(1/2+\ep)n$ (for some positive
constant $\ep$). Then

$$\sum_{a \in A, a < p/2} a \le p  + \frac{1}{\ep} (n+1) $$

\noindent and

$$ \sum_{a \in A, a > p/2} \|a \|  \le (\frac{1}{\ep}+1)n. $$
\end{lemma}

\begin{proof} (Proof of Lemma \ref{lemma:main2-1}.)  Following the proof of  Lemma \ref{lemma:simple1}, with
$l=|A''\cap I_n|$ and $k=|A''\setminus I_n|$, we have

$$(l+n+1)(n-l)>(2n+k)k.$$

\noindent On the other hand, $n\ge k+l=|A''|=|A|-O(n/\log^2 n),$
thus $n-l=k+n-|A|+O(n/\log^2 n)=(1-\lambda+o(1))n+k$ and
$n+l\le(1+\lambda)n-k.$ Putting all these together with the fact
that $\lambda$ is quite close to 1, we can conclude that  that $k
< .1n$. It follows (rather generously) that  $l=\lambda
n-k-O(n/\log^2 n)>.8n$.

The above shows that most of the elements of $A$ belong to $I_n$, as

$$|A_1|=|A\cap I_n|\ge |A''\cap I_n|> .8n.$$

\noindent Split $A_1$ into two sets, $A_1'$ and $A_1'':=A_1\setminus
A_1'$, where $A_1'$ contains all elements $a$ of $A_1$ such that
$n+1-a$ also belongs to $A_1$. Recall that $A_1$ has at least
$\lfloor n/2 \rfloor-t$ pairs $(a_i,b_i)$ satisfying $a_i+b_i=n+1$.
This guarantees that $|A_1'|\ge 2(\lfloor n/2 \rfloor-t)\ge.6 n$. On
the other hand, $A_1'$ is a subset of the core of $A$. The proof is
complete.
 \end{proof}

\begin{proof} (Proof of Lemma \ref{lemma:main2-2})
Abusing the notation slightly, we use $A_1'$ to denote the core of
$A$.  We have  $|A_1'|\ge (1/2+\epsilon)n$.

\begin{lemma}\label{lemma:simple8}
Any $l\in [n(1/\epsilon+1),n(1/\epsilon+1) + n]$ can be written as a
sum of $2(1/\epsilon+1)$ distinct elements of $A_1'$.
\end{lemma}

\begin{proof} First notice that for any $m$ belongs to
$I_{\epsilon}=[(1-\epsilon)n,(1+\epsilon)n]$, the number of pairs
$(a,b)\in {A_1'}^2$ satisfying $a<b$ and $a+b=m$ is at least
$\epsilon n/2$. Next, observe that any $k$, $k\in [0,n]$, is a sum
of $1/\epsilon +1$ integers (not necessarily distinct) from
$[0,\epsilon n]$. Consider $l$ from
$[n(1/\epsilon+1),n(1/\epsilon+1) + n]$; we can represent
$l-n(1/\epsilon+1)$ as a sum $a_1+\cdots+a_{1/\epsilon+1}$ where
$0\le a_1,\dots, a_{1/\epsilon+1}\le \epsilon n$. Thus $l$ can be
written as a sum of $1/\epsilon +1$ elements (not necessarily
distinct) of $I_{\epsilon}$, as
$l=(n+a_1)+\cdots+(n+a_{1/\epsilon+1}).$ Now we represent each
summand in the above representation of $l$ by two elements of
$A_1'$. By the first observation, the numbers of pairs are much
larger than the number of summands, we can manage so that all
elements of pairs are different.
\end{proof}

\noindent Recall that $A_1'$ consists of pairs $(a_i',b_i')$ where
$a_i'+b_i'=n+1$, so

$$\sum_{a'\in A_1'}a'=(n+1)|A_1'|/2.$$

\begin{lemma}\label{lemma:simple9}
$I':=[n(1/\epsilon+1), \sum_{a'\in A_1'}a'-(n+1)/\epsilon]\subset
S_{A_1'}.$
\end{lemma}

\begin{proof} Lemma \ref{lemma:simple8} implies that for each  $x \in
[n(1/\epsilon+1),n(1/\epsilon+1) + n]$ there exist distinct
elements $a_1',\dots, a_{2(1/\epsilon+1)}'\in A_1'$ such that
$x=\sum_{i=1}^{2(1/\epsilon+1)}a_i'$. We discard all $a_i'$ and
$(n+1)-a_i'$ from $A_1'$. Thus there remain exactly
$|A_1'|/2-2(1/\epsilon+1)$ different pairs $(a_i'',b_i'')$ where
$a_i''+b_i''=n+1$. The sums of these pairs represent all numbers
of the form $k(n+1)$ for any $0\le k \le
|A_1'|/2-2(1/\epsilon+1)$. We thus obtained a representation of
$x+k(n+1)$ as a sum of different elements of $A_1'$, in other word
$x+k(n+1)\in S_{A_1'}$. As $x$ varies in
$[n(1/\epsilon+1),n(1/\epsilon+1) + n]$ and $k$ varies in
$[0,|A_1'|/2-2(1/\epsilon+1)]$, the proof is completed.
\end{proof}

\noindent Let $A_2=A\setminus A_1$ and set $A_2':= A_2 \cap
[0,(p-1)/2]$ and $A_2'' = A_2 \backslash A_2'$. We are going to view
$A_2''$ as a subset of $[-(p-1)/2,-1]$.

\noindent We will now invoke Lemma \ref{lemma:simple4} several times
to conclude Lemma \ref{lemma:main2-2}. First, it is trivial that the
length of $I'$ is much larger than $n$, whilst elements of $A_1$ are
positive integers bounded by $n$. Thus, Lemma \ref{lemma:simple4}
implies that

$$I'':=[n(1/\epsilon+1),\sum_{a\in A_1}a-(n+1)/\epsilon ]\subset
S_{A_1}.$$

Note that the length of $I''$ is greater than $(p-1)/2$. Indeed $n
\approx (2p)^{1/2}$ and

$$|I''|=\sum_{a\in A_1}a-(n+1)/\epsilon
-n(1/\epsilon+1) \ge \sum_{a\in A_1'}a- O(n)$$
$$\ge (1/2+\epsilon)n(n+1)/2-O(n)>(p-1)/2.$$

Again, Lemma \ref{lemma:simple4} (applied to $I''$)  yields that

$$[n(1/\epsilon+1),\sum_{a\in A_1\cup A_2'}a-(n+1)/\epsilon ]\subset
S_{A_1\cup A_2'}$$

and

$$[\sum_{a\in A_2''}a+n(1/\epsilon+1),\sum_{a\in A_1}a-(n+1)/\epsilon ]\subset
S_{A_1\cup A_2''}.$$

The union of these two long intervals belongs to $S_A$

$$[\sum_{a\in A_2''}a+n(1/\epsilon+1),\sum_{a\in A_1\cup A_2'}a-(n+1)/\epsilon ]\subset S_A.$$

On the other hand, $0\notin S_A$ implies

$$\sum_{a\in A_2''}a+n(1/\epsilon+1)>0$$

and

$$\sum_{a\in A_1\cup A_2'}a-(n+1)/\epsilon < p.$$

\noindent The proof of Lemma \ref{lemma:main2-2} is completed.
\end{proof}

\section {Sketch of the proof of Theorem \ref{theorem:main3}}

Assume that $A$ is incomplete and $|A|=\lambda p^{1/2}$ with some
$2\ge \lambda\ge 1.99.$ Furthermore, assume that the element $b$ in
Theorem \ref{lemma:main2} is one. We are going to view $\BZ_p$ as
$[-(p-1)/2,(p-1)/2]$.

\noindent To make the proof simple, we made some new invention:
$n=\lfloor p^{1/2}\rfloor$, $A_1:=A\cap [-n,n], A_1':=A\cap [0,n],
A_1'':=A\cap [-n,-1], A_2':=A\cap [n+1,(p-1)/2], A_2'':=A\cap
[-(p-1)/2,-(n+1)], t_1':=|A_1'|, t_1'':=|A_1''|,
t_1:=|A_1|=t_1'+t_1''$.

\noindent Notice that $|A''|$ (in Theorem \ref{lemma:main2}) is
sufficiently close to the upper bound. The following holds.
\begin{lemma}\label{lemma:main3-1} Most of the elements of $A''(A)$
belong to $[-n,n]$;

\begin{itemize}
\item both $t_1'$ and $t_1''$ are larger than
$(1/2+\epsilon)n,$
\item $t_1$ is larger than $(2^{1/2}+\epsilon)n$
\end{itemize}

with some positive constant $\epsilon$.
\end{lemma}

\noindent As a consequent, both $S_{A\cap [-n,-1]}$ and $S_{A\cap
[1,n]}$ contain long intervals thanks to the following Lemma, which
is a direct application of Lemma \ref{lemma:simple8} and argument
provided in Lemma \ref{lemma:simple2}.

\begin{lemma}\label{lemma:main3-2}
If $X$ is a subset of $[1,n]$ with size at least $(1/2+\epsilon)n$.
Then
$$[(n+1)(1/\epsilon+1),(n+1)(n/2-t-c_{\epsilon})]\subset S_X$$
where $t=n-|X|$ and $c_{\epsilon}$ depends only on $\epsilon$.
\end{lemma}

\noindent Now we can invoke Lemma \ref{lemma:simple4} several times
to conclude Theorem \ref{theorem:main3}.

Lemma \ref{lemma:main3-2} implies

$$I':=[(n+1)(1/\epsilon+1),(n+1)(n/2-t_1'-c_{\epsilon})]\subset
S_{A_1'}.$$

and

$$I'':=[-(n+1)(n/2-t_1''-c_{\epsilon}),-(n+1)(1/\epsilon+1)]\subset
S_{A_1''}.$$

Lemma \ref{lemma:simple4} (applied to $I'$ and $A_1''$; $I''$ and
$A_1'$ respectively) yields

$$[\sum_{a_1''\in
A_1''}a_1''+(n+1)(1/\epsilon+1),(n+1)(n/2-t_1'-
c_{\epsilon})]\subset S_{A_1}$$

and

$$[-(n+1)(n/2-t_1''-c_{\epsilon}),\sum_{a_1'\in
A_1'}a_1'-(n+1)(1/\epsilon+1)]\subset S_{A_1}.$$

which gives

$$I:=[\sum_{a_1''\in A_1''}a_1''+(n+1)(1/\epsilon+1),\sum_{a_1'\in
A_1'}a_1'-(n+1)(1/\epsilon+1)]\subset S_{A_1}.$$

\noindent Note that the length of $I$ is greater than $(p-1)/2$.
Again, Lemma \ref{lemma:simple4} (applied to $I$ and $A_2'$, $I$ and
$A_2''$ respectively) implies

$$[\sum_{a''\in A_1''\cup A_2''}a''+(n+1)(1/\epsilon+1),\sum_{a_1'\in
A_1'}a_1'-(n+1)(1/\epsilon+1)] \subset S_A$$

and

$$[\sum_{a_1''\in A_1''}a_1''+(n+1)(1/\epsilon+1), \sum_{a'\in A_1'\cup
A_2'}a'-(n+1)(1/\epsilon+1)] \subset S_A.$$

The union of these two intervals belongs to $S_A$,

$$[\sum_{a''\in A_1''\cup A_2''}a''+(n+1)(1/\epsilon+1),\sum_{a'\in
A_1'\cup A_2'}a'-(n+1)(1/\epsilon+1)]\subset S_{A}.$$

On the other hand, $S_A\neq \BZ_p$ implies

$$\sum_{a'\in
A_1'\cup A_2'}a'-\sum_{a''\in A_1''\cup
A_2''}a''-2(n+1)(1/\epsilon+1)<p.$$

In other word

$$\sum_{a \in A} \| a\| \le p +O(p^{1/2}).$$

\end{document}